\documentclass{amsart}
\usepackage{graphicx}
\usepackage{amssymb}
\newtheorem{thm}{Theorem}[section]

\theoremstyle{definition}

\theoremstyle{question}
\newtheorem{que}[thm]{Question}
\theoremstyle{Conjecture}
\newtheorem{con}[thm]{Conjecture}

\numberwithin{equation}{section}

\begin{document}

\title[AN AFFIRMATIVE ANSWER TO a conjecture related to the solvability of groups ]{An affirmative answer to a conjecture related to  the solvability of  groups}%
\author{M. Zarrin}%

\address{Department of Mathematics, University of Kurdistan, P.O. Box: 416, Sanandaj, Iran}%
 \email{M.zarrin@uok.ac.ir, zarrin@ipm.ir}
\begin{abstract}
 In this paper, we show that each finite group $G$ containing at most $p^2$ Sylow $p$-subgroups for each odd prime number $p$,  is a solvable group. In fact, we give a positive answer to the conjecture in \cite{Rob}. 
\end{abstract}
\maketitle

\section{\textbf{ Introduction}}

we know that, according to Feit-Thompson Theorem,  every group with an odd order is solvable. As a consequence of this Theorem, one can say that every finite group $G$ that has a normal Sylow $2$-subgroup, i.e., with $v_2(G) = 1$, is a solvable group, where $v_p(G)$ is the number of Sylow p-subgroups of $G$.  Also, if $v_p(G)=1$ for each prime $p$, then $G$ is nilpotent and reciprocally. This result show that the number of Sylow $p$-subgroups for a prime p is restricted arithmetically by the properties of a  group $G$.  Most recently, the author in \cite{Rob}, proved that if $v_p(G)\leq p^2-p+1$ for each prime p, then $G$ is solvable. Here, first we show that it is not necessary to consider the  amount of $v_p(G)$ for each prime number $p$ and also improve the upper bound to $p^2$.  In fact, we give  a substantial generalization as follows:

\begin{thm}
Every finite group $G$ containing at most $4$ Sylow $2$-subgroups,  is a solvable group.
\end{thm} 

Also, the author,  raised the following conjecture.
\begin{con} 
Let G be a finite group. If  $v_p(G)\leq p^2-p+1$ for each odd prime number $p$, then $G$ is solvable.
\end{con} 
Finally,  we  give the positive  answer to this conjecture  and improve it as follows:  

\begin{thm}
Every finite group $G$ containing at most $p^2$ Sylow $p$-subgroups for each odd prime number $p$,  is solvable. 
\end{thm}

\section{\textbf{The Proofs}}

$\mathbf{Proof~~ of ~~Theorem ~~1.1}.$

Suppose, on the contrary, that there exists a non-solvable finite group $G$ of the least possible order with $v_2(G)\leq 4$. In this case $G$ should be a simple group. Otherwise, if there exists a non-trivial proper normal subgroup $M$ of $G$, then as $v_2(M)\leq v_2(G)\leq 4$ and $ v_2(G/M)\leq v_2(G)\leq 4$, both $M$ and $G/M$ are soluble (note that if $N$ or $G/N$ is  group with odd order, then by Feit-Thompson Theorem they are solvable). It follows that $G$ is solvable, which is a contradiction. Therefore $G$ is a minimal simple group with $v_2(G)\leq 4$. By Thompson's classification of minimal simple groups \cite{Tho}, $G$ is isomorphic to one of the following simple groups: $A_5$  the alternating group of degree $5$, $L_2(2^p )$, where $p$ is an odd prime,
$L_2(3^p)$, where $p$ is an odd prime,
$L_2(p)$, where $5 < p$ is prime and $p \equiv 2 (mod~~5)$,
$L_3(3)$, and $^2B_2(q)$ where $q=2^{2m+1}\geq 8$.\\
Now we show that in each case we obtain a contradiction. This completes the proof.
Clearly $v_2(A_5)=5$ and $v_2(L_3(3))=351$, a contradiction.

If $G$ is isomorphic to $L_2(2^p)$, then by Case 2 of the proof of Proposition 2.4 of \cite{Shi},  we get  that $n_2(G)=2^p+1\geq 9$, a contradiction.

If $G$ is isomorphic to $L_2(3^p)$, then one can again imply,  from  Proposition 2.4 of \cite{Shi}, that  $5 < v_2(G)=3^{2p}-1$ or $(3^{3p}- 3^{p})/24$, a contradiction.

If $G$ is isomorphic to $L_2(p)$, where $5 < p$ is prime and $p = 2 (mod~~5)$, then by an argument similar to $L_2(3^p)$ we obtain that  $5 < v_2(G)=p^{2}-1$ or $(p^{3}- p)/24$, a contradiction.

If $^2B_2(q)$,  $q=2^{p}$ and $p$ an odd prime, then  by Theorem 3.10 (and its proof) of Chapter XI of \cite{Hup}, we have $|G| = (q-1)(q^{2})(q^{2} +1)$ and $v_2(G) = q^{2} + 1> 65$, a contradiction. \\

We note that the bound 4 in Theorem 1.1 is the best possible, as $v_2(A_5)=5$.

Now by similar argument we prove Theorem 1.3.\\\\
$\mathbf{Proof~~ of ~~Theorem ~~1.3}.$

Suppose, on the contrary, that there exists a non-solvable finite group $G$ of the least possible order with $v_p(G)\leq p^2$ for all its prime odd divisors. In the sequel,  by an argument similar to the proof of Theorem 1.1, to prove it is enough to consider the following groups (note that $v_3(A_5)=10$, $v_3(L_3(3))=52$): \\
If $G$ is isomorphic to $L_2(q)$ with $q=2^p$, then we consider  an odd prime divisor of $|G|$, like $r$.  Then it is easy to see that  $r$ divides  either  $q+1$ or $q-1$. Now if $R$ is a Sylow $r$-subgroup of $G$, then $R$ is cyclic such that $N_G(R)=D_{q-1}$ or $D_{q+1}$, where $D_m$ is the dihedral group of order $m$. Therefore, the number of Sylow $r$-subgroups is $q(q + 1)/2$ or $q(q - 1)/2$ and so  $n_r(G)> r^2$, a contradiction.

If $G$ is isomorphic to $L_2(q)$ with $q=3^p$ and $p$ is an odd prime, then it is easy to see that $$v_3(L_2(q))=v_3(SL(2,q)/Z(SL(2,q)))=v_3(SL(2,q)),$$where $Z(SL(2,q))$ is the center of the group $SL(2,q)$. Assume that $R \in Syl_G(3)$, then 
 $N_G(R)$ is the set of upper triangular matrices with determinant 1. Therefore,  the order of the normalizer $N_G(R)$ is $q(q-1)$. Thus $v_3(G)=q(q^2 -1)/q(q-1)=q +1 > 3^2$, a contradiction.
 
If $G$ is isomorphic to $L_2(p)$, where $5 < p$ is prime and $p\equiv 2 (mod~~5)$, then by an argument similar to $L_2(2^p)$ we obtain that  $v_r(G)=q(q+1)/2$ or $q(q-1)/2$, where $r\neq 2$, and so  $n_r(G)> r^2$, a contradiction.

If $G=^2B_2(q)$, where  $q=2^{2m+1}\geq 8$, then it is well-known that the Suzuki group $^2B_2(q)$ contains a maximal subgroup like $T$ of order $4(q-r+1)$, where $r=2^{m+1}$ and also $T$ has a normal cyclic subgroup $C$ in which $|C|=q-r+1$. Moreover, $C$ includes a Sylow $5$-subgroup like $P$ (note that as $q^2+1\equiv 0 (mod ~5)$ so  $5$ is a prime divisor of $|G|$). From this one can follow that $T \leq N_G(P)$ and so $T=N_G(P)$, as $T$ is maximal. Since $|G| = (q-1)(q^{2})(q^{2} +1)$,  the number of conjugates of $P$ in $G$ is $$v_5(G)=|G|/|T|=(q-1)(q^{2})(q^{2} +1)/4(q-r+1)> 25.$$ Thus in each case we obtain a contradiction. This completes the proof.\\

Finally, it is well-known  that the only nonabelian simple finite groups in which its order is not divisible by 3 are the Suzuki  groups. From this one can show, by induction on the order, that:  if  $H$ is a group such that $v_3(H) = 1$ and has no composition factor isomorphic to $^2B_2(q)$, then $H$ is a solvable group. As a result,  one can find out that some of the odd prime numbers (for instance, 3) have stronger influence on the solvability of groups. In fact, most probably  for the solvability of finite groups in terms of the number of Sylow $p$ subgroups, we do not need to consider all odd prime numbers. Therefore, it might seem reasonable to pose the following question:

\begin{que} 
 What is the smallest positive integer $n$ such that whenever there exist a finite group $G$ satisfying $v_{p_i}(G )\leq p_i^2$ where $p_i$ is odd number and $i\in \{1,\dots, n\}$, which guarantees  the  solvability of $G$?
\end{que}

\end{document}